\documentclass[preprint,12pt]{elsarticle}
\usepackage{geometry}
\usepackage{orcidlink}
\usepackage{amssymb, amsmath}
\usepackage{amsthm}
\usepackage[utf8]{inputenc}

\usepackage{hyperref}
\usepackage{enumitem}
\hypersetup{
  colorlinks=true,
}
\hypersetup{
    linkcolor=magenta,
    filecolor=magenta,      
    urlcolor=cyan,
    citecolor=magenta,
    }

\theoremstyle{remark}

\theoremstyle{definition}

\newtheorem{definition}{Definition}

\journal{xxxxxxxxxxx}

\newcommand{\cinf}{$\mathcal{C}^{\infty}$}
\newcommand{\contract}{\,\lrcorner\,}

\usepackage{todonotes}
\usepackage{xcolor}

\begin{document}
\hypersetup{
    linkcolor=magenta,
    filecolor=magenta,      
    urlcolor=cyan,
    citecolor=magenta,
    }
\begin{frontmatter}



\title{A \texorpdfstring{\cinf}{Cinf}-structures approach to finding traveling wave solutions of the gKdV equation}


\author[pan]{Antonio J. Pan-Collantes\orcidlink{0000-0002-0063-0710}}

\affiliation[pan]{
organization={Departamento de Matematicas,
Universidad de Cadiz - UCA},
            addressline={Facultad de Ciencias, Campus Universitario de Puerto Real s/n}, 
            city={Puerto Real},
            postcode={11510}, 
            state={Cadiz},
            country={Spain}}

\begin{abstract}
A novel geometric method is applied to the problem of describing traveling wave solutions of the generalized Korteweg–de Vries (gKdV) equation in the form
$$
u_t + u_{xxx} + a(u)u_x = 0,
$$
where $a(u)$ is a smooth function characterizing the nonlinearity. Using the traveling wave ansatz, the gKdV equation reduces to an ordinary differential equation (ODE), which we analyze via the $\mathcal{C}^\infty$-structure-based method, a geometric framework involving sequences of involutive distributions and Pfaffian equations. Starting with the symmetry $\partial_z$, we construct a $\mathcal{C}^\infty$-structure for the ODE and apply the stepwise integration algorithm to obtain an implicit general solution. Then we derive explicit solutions for specific forms of $a(u)$, including the modified KdV and Schamel--KdV equations, as well as power-law nonlinearities.

\end{abstract}

\begin{keyword}
\cinf-structure \sep traveling wave solutions \sep Korteweg--de Vries \sep integrable systems

\MSC[2020] 35Q53 \sep 35C07 \sep 37K10 \sep 58A17
\end{keyword}

\end{frontmatter}


\section{Introduction}

The Korteweg--de Vries (KdV) equation is a significant contribution to mathematics and science. It is a nonlinear partial differential equation that describes the evolution of waves on shallow water surfaces. Its origins are traced back to 1834, when the naval engineer J. S. Russell observed a wave in the Union Canal in Scotland that maintained its shape and velocity while traveling through the channel \cite{russell1844report}. 

Dutch mathematicians D. J. Korteweg and his student G. de Vries derived a nonlinear partial differential equation that could be used to model this kind of behavior \cite{korteweg1895xli}. After a change of variable, their equation can be rewritten as
$$
u_t+u_{xxx} + 6uu_x = 0.
$$

This equation, nowadays known as the Korteweg--de Vries equation (KdV), describes the evolution of a one-dimensional wave in a dispersive medium. It has led to the discovery of new mathematical structures and techniques, establishing a foundational moment in the study of integrable systems. Because of its well-understood soliton solutions, conservation laws, and integrable structure, the KdV equation is also frequently used as a benchmark to test new analytical and numerical integration methods.

Over the years, the KdV equation has been generalized in numerous ways to incorporate more complex phenomena and to enable the study of higher-dimensional and more intricate wave structures. In this paper, we will work with a generalized KdV (gKdV) equation with the form
\begin{equation}\label{c6EqPDEKdV}
u_t+u_{xxx} + a(u)u_x = 0,
\end{equation}
where $a\in \mathcal{C}^{\infty}(\mathbb R)$, containing the original KdV equation as a particular case \cite{kdvTsutsumi,kdvKenig,kdvChrist}. This family of equations has attracted significant attention due to its rich mathematical structure and diverse applications in different areas, for example: 
\begin{itemize}
	\item For $a(u)=u^2$, equation \eqref{c6EqPDEKdV} is known as the modified Korteweg--de Vries equation, and it is used to describe the behavior of warm plasma with negative ions \cite{kalita1990modified} and dusty plasmas \cite{emami2011solitons}. 

	\item For $a(u)=\alpha \sqrt{u}$, equation \eqref{c6EqPDEKdV} becomes the Schamel equation, which describes the propagation of ion-acoustic waves in a cold-ion plasma, in which specific electrons are trapped and do not exhibit isothermal behavior during the wave's passage \cite{tariq2022new}. 
	
	\item \sloppy For $a(u)=\alpha \sqrt{u}+\beta u$, with $\alpha, \beta \in \mathbb R$, equation \eqref{c6EqPDEKdV} corresponds to the Schamel--Korteweg--de Vries equation \cite{giresunlu2017exact}, and for ${a(u)=1+\alpha \sqrt{u}+\beta \ln(|u|)}$ the logarithmic Schamel equation. They both model different trapping scenarios for the electrons. 
	
	\item For $a(u)=2\alpha u-\beta u^2$, $\alpha, \beta \in \mathbb R$, equation \eqref{c6EqPDEKdV} is called Gardner equation, and it has applications in hydrodynamics, plasma physics and quantum field theory \cite{wazwaz2007new,betchewe2013new}.
	
	\item Finally, for $a(u)=u\ln(|u|)$ equation \eqref{c6EqPDEKdV} is called the Logarithmic KdV equation, which describes solitary waves in harmonic chains with Hertzian interaction force \cite{inc2018investigation}.
	
\end{itemize}

In this work, we address the problem of finding exact solutions to equation \eqref{c6EqPDEKdV}, and more specifically traveling wave solutions. These are solutions whose shape remains unchanged while they propagate at a constant speed. To this end, we reduce equation \eqref{c6EqPDEKdV} to an ordinary differential equation (ODE), and then apply the \cinf-structure-based method of integration, which is a geometric method designed to find the integral manifolds of involutive distributions \cite{pancinf-sym,pancinf-struct,pan23integration}.

It is important to clarify the primary contribution of this work. We acknowledge that for several of the specific cases examined in the paper, the resulting ODEs may be solved through more direct integration techniques. However, the goal here is not to present a computational shortcut for these well-known equations. Instead, our aim is to showcase the efficacy of the \cinf-structure formalism as a systematic, reproducible pipeline, with potential to address other nonlinear differential equations where direct solutions may be not immediately apparent.

The structure of the paper is as follows. In Section \ref{cinfbasesmethod}, we briefly review the \cinf-structure-based method of integration. In Section \ref{travsolutions}, we apply this method to the gKdV equation, obtaining an implicit general solution. In Section \ref{particularcases}, we analyze some particular cases of the function $a(u)$, obtaining explicit traveling wave solutions to the modified KdV equation, the Schamel--KdV equation and the logarithmic Schamel equation.

%%%%%%%%%%%%%%%%%%%%%%%%%%%%%%%%%%%%%%%%%%%%%%
\section{\texorpdfstring{\cinf}{cinf}-structure-based method of integration}\label{cinfbasesmethod}
%%%%%%%%%%%%%%%%%%%%%%%%%%%%%%%%%%%%%%%%%%%%%%

The \cinf-structure-based method of integration is a novel approach to find integral manifolds of involutive distributions. It has been successfully applied to integrate high-order ODEs that do not have enough Lie point symmetries for the classical methods, and also to dynamical systems and partial differential equations \cite{pancinf-sym,pancinf-struct,pan23integration,PANCOLLANTES2025116091}.

Given an open subset $U\subset \mathbb{R}^n$, we define a distribution on $U$ of rank $r$ as a submodule of the $\mathcal C^{\infty}(U)$-module of vector fields $\mathfrak{X}(U)$ generated by a set of $r$ pointwise linearly independent vector fields $\{Z_1,\ldots,Z_r\}$. We will denote this distribution by $\mathcal{S}(\{Z_1,\ldots,Z_r\})$, and it will be said to be involutive if 
$$
[Z_i,Z_j]\in \mathcal{S}(\{Z_1,\ldots,Z_r\})
$$ 
for all $i,j=1,\ldots,r$, where $[-,-]$ denotes the Lie bracket of vector fields \cite{warner,lee2013smooth}.

In order to apply the method to the particular case of differential equations, consider, for $m\geq 2$, the jet bundle $J^{m-1}(\mathbb R, \mathbb R)$ with coordinates $(z,y,y_1,\ldots,y_{m-1})$, where $y_i=\frac{d^i y}{dz^i}$ for $i=1,\ldots,m-1$ (see \cite{saunders1989geometry} for details). Define the standard volume form on $J^{m-1}(\mathbb R, \mathbb R)$ as
$$
\boldsymbol{\Omega} = dz \wedge dy \wedge d y_1 \wedge \cdots \wedge dy_{m-1}.
$$
In this geometric setting, an $m$th-order ODE is an expression of the form
	\begin{equation}\label{ode_m}
		y_m=\phi(z,y,y_1,\ldots,y_{m-1}),
	\end{equation}
where $\phi$ is a smooth function defined on an open subset $U\subseteq J^{m-1}(\mathbb R, \mathbb R)$. We can associate to equation \eqref{ode_m} the (trivially) involutive distribution $\mathcal{S}(\{Z\})$, generated by the vector field
$$
Z=\partial_z+y_1\partial_y+\cdots+\phi \partial_{y_{m-1}},
$$ 
defined on $U$.

The integral manifolds of this distribution correspond to the prolongation of solutions of the ODE \eqref{ode_m} to the jet space $J^{m-1}(\mathbb R, \mathbb R)$ \cite{saunders1989geometry}. To look for these integral manifolds by using the \cinf-structure-based method we need to introduce the notion of \cinf-structure \cite{pancinf-sym,pancinf-struct}, which is a generalization of the concept of solvable structure \cite{basarab,hartl1994solvable}.

\begin{definition}\label{Qsolvable}
An ordered collection of vector fields $( X_1,\ldots,X_{m})$ is a \cinf-structure for equation \eqref{ode_m} if the distribution
	$$
\mathcal{S}(\{Z,X_1,\ldots,X_i\})
	$$
has constant rank $i+1$ and is involutive, for $1\leq i\leq m$.
\end{definition}

The main result concerning  \cinf-structures is that they can be used to integrate the ODE \eqref{ode_m} by solving $m$ completely integrable Pfaffian equations (see \cite[Theorem 4.1]{pancinf-sym}). For the reader’s convenience, we present an outline of the integration procedure below:

\begin{enumerate}
	\item We define the 1-forms 
	\begin{equation}\label{formas}
		\omega_i={X_{m}\,\lrcorner\,\ldots\,\lrcorner\,\widehat{X_i}\,\lrcorner\,\ldots\,\lrcorner\,X_{1}\,\lrcorner\,Z\,\lrcorner\,\boldsymbol{\Omega}}, \quad 1\leq i\leq m,
		\end{equation}
		where $\widehat{X_i}$ indicates omission of $X_i$ and $\lrcorner$ denotes interior product.

		\item The Pfaffian equation $\omega_{m}\equiv 0$ is completely integrable (see \cite{pancinf-sym,pancinf-struct,pan23integration}). A solution $I_m=I_m(z,y,y_1,\ldots,y_{m-1})$ can be found by solving the system of linear homogeneous first-order PDEs arising from the condition $dI_m\wedge \omega_m=0$. 
		
		\item \sloppy Once the function $I_m$ is found, we consider the level sets $\Sigma_{(C_m)}$ implicitly defined by ${I_m(z,y,y_1,\ldots,y_{m-1})=C_m}$, where $C_m\in \mathbb R$. And we look for a local parametrization $\iota_m$ of $\Sigma_{(C_m)}$, with domain an open subset $U_m\subseteq \mathbb R^{m}$.
		
		\item We compute the pullbacks
		$$
		\omega_i|_{\Sigma_{(C_m)}}:=\iota_m^*(\omega_i), \quad 1\leq i\leq m.
		$$
		Observe that $\omega_m|_{\Sigma_{(C_m)}}= 0$.
		
		\item We repeat the process for the Pfaffian equation $\omega_{m-1}|_{\Sigma_{(C_m)}}\equiv 0$, which is now completely integrable, and it is defined on a space of one dimension less, that is, on $U_{m}$ instead of $J^{m-1}(\mathbb R, \mathbb R)$. The procedure gives rise to new 1-forms
		$$
		\omega_{i}|_{\Sigma_{(C_{m-1},C_m)}}:=\iota_{m-1}^*(\omega_{i}|_{\Sigma_{(C_m)}}), \quad 1\leq i\leq m,
		$$ 
		defined on certain open set $U_{m-1} \subseteq \mathbb R^{m-1}$.

		Now, both $\omega_{m}|_{\Sigma_{(C_{m-1},C_m)}}$ and $\omega_{m-1}|_{\Sigma_{(C_{m-1},C_m)}}$ are null 1-forms, so we proceed with the Pfaffian equation $\omega_{m-2}|_{\Sigma_{(C_{m-1},C_m)}}\equiv 0$, which is now completely integrable.

		\item \sloppy In each iteration, for $m-1 \geq k\geq 1$, the completely integrable Pfaffian equation 
		$$
		\omega_k|_{\Sigma_{(C_{k+1},\ldots,C_m)}}\equiv 0
		$$
		is solved, and a solution 
		$$
		I_k=I_k(z,y,y_1,\ldots,y_{k-1};C_{k+1},\ldots,C_{m})
		$$ 
		is obtained. Then, a local parametrization $\iota_k$, defined on an open subset $U_k\subseteq \mathbb R^k$, of the level set $\Sigma_{(C_k,\ldots,C_m)}$, $C_k\in \mathbb R$, is obtained. The process terminates when a local parametrization $\iota_1$ of the level set $\Sigma_{(C_1,\ldots,C_m)}$ is obtained.
		
		\item Finally, the composition of local parametrizations, $\iota_m\circ\cdots\circ\iota_1$, yields a parametric expression for the integral curves of $Z$, which are the prolongations, in the jet space $J^{m-1}(\mathbb R, \mathbb R)$, of the explicit solutions to the ODE \eqref{ode_m}. Alternatively, the expression 
		$$
		I_1(z,y;C_2,\ldots,C_m)=C_1,\quad C_1\in \mathbb R,
		$$
		provides an implicit general solution to the ODE \eqref{ode_m}.
\end{enumerate}

Readers interested in a more detailed explanation of the \cinf-structure integration process, and in reviewing illustrative examples, are referred to \cite{pancinf-sym,pancinf-struct,pan23integration}.

\section{Traveling wave solutions to the gKdV equation}\label{travsolutions}
%%%%%%%%%%%%%%%%%%%%%%%%%%%%%%%%%%%%%%%%%%%%%%%%%%%%%%%%%%%%%%
We now consider the problem of finding traveling wave solutions of equation \eqref{c6EqPDEKdV}. To this end, we employ the standard ansatz consisting of the transformation
\begin{equation}\label{c6EqTrav2D}
\begin{array}{l}
z = x - c t, \
y(z) = u(x,t),
\end{array}
\end{equation}
which reduces the PDE \eqref{c6EqPDEKdV} to the third-order ODE
\begin{equation}\label{c6EqODEKdV}
-c y_1 + y_3 + a(y)y_1 = 0.
\end{equation}

% Now we face the problem of looking for the traveling wave solutions to equation \eqref{c6EqPDEKdV}. To seek traveling wave solutions, we employ the standard ansatz, which consists of the transformation
% \begin{equation}\label{c6EqTrav2D}
% 	\begin{array}{l}
% 		z=x-c t, \\	
% 		y(z)=u(x,t),	
% 	\end{array}
% \end{equation}
% which allows us to reduce the PDE \eqref{c6EqPDEKdV} to the third-order ODE
% \begin{equation}\label{c6EqODEKdV}
% 	-c y_1+y_3+a(y)y_1=0.
% \end{equation}

The corresponding vector field in the jet space $J^2(\mathbb R, \mathbb R)$ is
$$
Z=\partial_z+y_1\partial_y+y_2 \partial_{y_1}+(c y_1-a(y) y_1)\partial_{y_2}.
$$

In order to find exact solutions to \eqref{c6EqODEKdV} using the method proposed in Section \ref{cinfbasesmethod}, we first need to identify a suitable \cinf-structure. We begin by observing that the vector field $\partial_z$ is a Lie symmetry of \eqref{c6EqODEKdV}; indeed, $[\partial_z,Z]=0$, as can be easily verified. Therefore, $\mathcal{S}(\{Z,\partial_z\})$ is involutive, and we can use the symmetry $X_1=\partial_z$ as the first vector field in our \cinf-structure.

To find the second vector field $X_2$, we consider the following ansatz:
$$
X_2=\partial_{y_1}+\eta \partial_{y_2},
$$
where $\eta=\eta(z,y,y_1)$ is a smooth function to be determined. According to Definition \ref{Qsolvable}, the distribution $\mathcal{S}(\{Z,X_1,X_2\})$ must be involutive, which requires that the Lie brackets $[Z,X_2]$ and $[X_1,X_2]$ lie in $\mathcal{S}(\{Z,X_1,X_2\})$. We compute these Lie brackets and impose that they are linear combinations of $Z$, $X_1$ and $X_2$, obtaining the following determining equations for $\eta$:
\begin{equation}\label{determiningPDE}
	\begin{aligned}
	\eta^2 y_1+\eta_y y_1^2+\eta_{y_1} y_1 y_2-\eta y_2+\eta_z y_1=0,\\
	\eta_z=0.		
	\end{aligned}
\end{equation}
A particular solution $\eta=\dfrac{y_1}{y}$ for the PDE system \eqref{determiningPDE} is obtained, and then the second vector field will be
$$
X_2=\partial_{y_1}+\dfrac{y_1}{y}\partial_{y_2}.
$$

For the last vector field in the \cinf-structure, we may choose any vector field $X_3$ that is linearly independent of $Z, X_1$ and $X_2$, since the distribution $\mathcal S (\{Z,X_1,X_2,X_3\})$ is then trivially involutive. A natural choice is $X_3=\partial_{y_2}$. Thus, we have a \cinf-structure $( X_1,X_2,X_3)$ for equation \eqref{c6EqODEKdV} given by
\begin{equation}\label{c6EqCanZK2}
	\begin{aligned}
		X_1&=\partial_z,\\
		X_2&=\partial_{y_1}+\dfrac{y_1}{y}\partial_{y_2},\\
		X_3&=\partial_{y_2}.
	\end{aligned}
\end{equation}

Now, we proceed by using the method described in Section \ref{cinfbasesmethod} to integrate the ODE \eqref{c6EqODEKdV}. First, we calculate the 1-forms \eqref{formas}
\begin{equation}\label{kdv1formas}
	\begin{aligned}
		\omega_1&= X_3\contract X_2\contract Z\contract \boldsymbol{\Omega}=-y_1 dz + dy,\\
		\omega_2&= X_3\contract X_1\contract Z\contract \boldsymbol{\Omega}=-y_2 dy +y_1 dy_1, \\
		\omega_3&= X_2\contract X_1\contract Z\contract \boldsymbol{\Omega}=-y_1\left(\dfrac{y_2}{y}+a(y)-c\right) dy + \dfrac{y_1^2}{y}  dy_1 - y_1dy_2.
		\end{aligned}
\end{equation}

A solution to the completely integrable Pfaffian equation $\omega_3\equiv 0$ is the smooth function defined by
$$
I_3=yy_2 - \frac{1}{2}y_1^2 -H_1(y),
$$
where $H_1(y)$ is any smooth function satisfying 
\begin{equation}\label{laH1}
	H_1'(y)=y(c-a(y)).
\end{equation}

We define the level sets
$$
\Sigma_{(C_3)}=\left\{(z,y,y_1,y_2)\in \mathbb R^4 \;\middle|\; I_3(z,y,y_1,y_2)=C_3\right\},
$$
with $C_3\in \mathbb R$. The restriction of the 1-forms \eqref{kdv1formas} to these level sets are
\begin{equation}
	\label{kdv1formas2}
	\begin{aligned}
	&\omega_1|_{\Sigma_{(C_3)}}=\iota_3^*(\omega_1)=-y_1 dz + dy,\\
	&\omega_2|_{\Sigma_{(C_3)}}= \iota_3^*(\omega_2)=-\frac{ y_1^2 + 2H_1(y)  + 2C_3}{2y}dy +y_1 dy_1,\\
	&\omega_3|_{\Sigma_{(C_3)}}= \iota_3^*(\omega_3)=0,
\end{aligned}
\end{equation}
where $\iota_3$ is the local parametrization of $\Sigma_{(C_3)}$ given by
$$
\iota_3(z,y,y_1)=\left(z,y,y_1, \dfrac{y_1^2 +2H_1(y)+2C_3}{2y}\right), \quad y\neq 0.
$$

We follow by looking for a solution to the Pfaffian equation $\omega_2|_{\Sigma_{(C_3)}}\equiv 0$, and we obtain the expression
$$
I_2=\dfrac{y_1^2-2yH_2(y)+2 C_3}{y},  \quad y\neq 0,
$$
where $H_2$ is any smooth function satisfying 
\begin{equation}\label{laH2}
	H_2'(y)=\dfrac{H_1(y)}{y^2}.
\end{equation}

We restrict again to level sets, in this case:
$$
\Sigma_{(C_2, C_3)} = \left\{ (x, y, y_1) \in \mathbb{R}^3 \;\middle|\; I_2(z, y, y_1; C_3) = C_2 \text{ and } y \ne 0 \right\},\quad C_2\in \mathbb R,
$$
and pull back the 1-forms \eqref{kdv1formas2} with the local parametrization
\begin{equation}\label{laiota2}
	\iota_2(z,y;C_3)=\left(z,y,\sqrt{y (C_2 + 2 H_2(y)) - 2 C_3}\right).
\end{equation}
In this way, we obtain
\begin{equation}\label{ultima1forma}
	\begin{aligned}
	&\omega_1|_{\Sigma_{(C_2,C_3)}}=-\sqrt{y (C_2 + 2 H_2(y)) - 2 C_3} dz + dy.\\
	&\omega_2|_{\Sigma_{(C_2,C_3)}}=0,\\
	&\omega_3|_{\Sigma_{(C_2,C_3)}}=0,
\end{aligned}
\end{equation}

Finally, a solution for the Pfaffian equation $\omega_1|_{\Sigma_{(C_2,C_3)}}\equiv 0$ is
$$
I_1= z-H_3(y;C_2,C_3),
$$
where $H_3$ satisfies
\begin{equation}\label{laH3}
	H_3'(y;C_2,C_3)=\dfrac{1}{\sqrt{y (C_2 + 2 H_2(y)) - 2 C_3}}.	
\end{equation}

Observe that the expression 
\begin{equation}\label{laI3}
z-H_3(y;C_2,C_3)=C_1
\end{equation}
is an implicit general solution to the ODE \eqref{c6EqODEKdV}. Then the traveling wave solutions to the gKdV equation \eqref{c6EqPDEKdV} are given by
\begin{equation}\label{TWSgKdV}
u(x,t)=F(x-ct-C_1;C_2,C_3),	
\end{equation}
where $F$ is an inverse function of $H_3$.

This is a very general expression, which depends not only on the parameters $c,C_1,C_2,C_3$, but also crucially on the choice of the function $a(u)$ in the definition of equation~\eqref{c6EqPDEKdV}. The specific choice of $a(u)$ determines the explicit forms of $H_1,H_2,H_3$, and consequently for $F$. To illustrate how different forms of $a(u)$ give rise to distinct classes of solutions, we now turn to the analysis of several particular cases.

%%%%%%%%%%%%%%%%%%%%%%%%%%%%%%%%%%%%%%%%%%%%%%%%%%%%%%%%%%%%%%
\section{Particular cases of \texorpdfstring{$a(u)$}{a}}\label{particularcases}
%%%%%%%%%%%%%%%%%%%%%%%%%%%%%%%%%%%%%%%%%%%%%%%%%%%%%%%
%%%%%%%%%%%%%%%%%%%%%%%%%%%%%%%%%%%%%%%%%%%%%%%%%%%%%%%%%

We will focus on deriving explicit traveling wave solutions to \eqref{c6EqPDEKdV} for particular forms of the function $a(u)$, beginning with the classical case:

\subsection{Case \texorpdfstring{$a(u)=6u$}{a(u)=6u}.}
%%%%%%%%%%%%%%%%%%%%%%%%%%%%%%

In this case, the PDE \eqref{c6EqPDEKdV} is the usual KdV equation, and the final Pfaffian equation $\omega_1|_{\Sigma_{(C_2,C_3)}}\equiv 0$ given in \eqref{ultima1forma}	becomes
	\begin{equation}\label{lPfafKdV}
		-\sqrt{-2y^3+cy^2+C_2 y-2C_3} \,dz + dy \equiv 0.
	\end{equation}
	Observe that this Pfaffian equation can be solved by quadrature, by evaluating
	\begin{equation}\label{classicKdVprimitive}
		\int \frac{1}{\sqrt{-2y^3+cy^2+C_2 y-2C_3}} \,dy,
	\end{equation}
	but this calculation is not straightforward. So we assume $C_2=C_3=0$, and we find the solution to \eqref{lPfafKdV} defined by
	$$
	I_1= z+\dfrac{2}{\sqrt{c}}\arctan\left(\sqrt{\dfrac{c-2 y}{c}}\right)=C_1,
	$$
	for $c>0$ and $\dfrac{c-2 y}{c}>0$; and the solution
	$$
	I_1= z-\dfrac{2}{\sqrt{-c}}\arctan\left(\sqrt{\dfrac{2 y-c}{c}}\right)=C_1,
	$$
	for $c<0$ and $\dfrac{2 y-c}{c}>0$.
	
Then, for $c>0$, a 1-parameter family of solutions to \eqref{c6EqODEKdV} is
	$$
	y(z)=\dfrac{c}{2}\mbox{sech}^2\left(\dfrac{\sqrt{c}}{2}(C_1-z)\right).
	$$

	Respectively, for $c<0$, a 1-parameter family of solutions to \eqref{c6EqODEKdV} is
	$$
	y(z)=\dfrac{c}{2}\mbox{sec}^2\left(\dfrac{\sqrt{-c}}{2}(C_1-z)\right).
	$$
	
	By means of the transformation \eqref{c6EqTrav2D}, these two families give rise, respectively, to the following two families of traveling wave solutions to the KdV equation:
	$$
	u(x,t)=\dfrac{c}{2}\mbox{sech}^2\left(\dfrac{\sqrt{c}}{2}(-x+ct+C_1)\right),
	$$
	and
	$$
	u(x,t)=\dfrac{c}{2}\mbox{sec}^2\left(\dfrac{\sqrt{-c}}{2}(-x+ct+C_1)\right).
	$$

% 	\begin{remark}\label{lodelaZK}
% It is worth noting that the search for the primitive \eqref{classicKdVprimitive} can be further analyzed through the roots of the polynomial $P(y) = -2y^3 + cy^2 + C_2 y - 2C_3$. In this sense, a complete description of the traveling wave solutions to the KdV equation may be obtained by employing an approach similar to that in \cite{PANCOLLANTES2025116091}. This analysis is beyond the scope of the present paper.
% \end{remark}

\subsection{Case \texorpdfstring{$a(u)=u^2$}{a(u)=u2}.}
%%%%%%%%%%%%%%%%%%%%%%%%%%%%%%

The corresponding equation \eqref{c6EqPDEKdV} is the modified KdV equation, and the final Pfaffian equation \eqref{ultima1forma}, $\omega_1|_{\Sigma_{(C_2,C_3)}}\equiv 0$, becomes
	$$
	-\frac{1}{6}\sqrt{-6y^{4}+36cy^{2}+36C_{2}y-72C_{3}}dz+dy\equiv 0.
	$$
	This Pfaffian equation can be solved by quadrature, and evaluating
	\begin{equation}\label{primitivemodifKdV}
		\int \frac{-6}{\sqrt{-6y^{4}+36cy^{2}+36C_{2}y-72C_{3}}} dy
	\end{equation}
	we can obtain implicit descriptions of the solutions to \eqref{c6EqODEKdV} based on the incomplete elliptic integral of first kind \cite{handbookfunctions}.

On the other hand, we can easily obtain an explicit expression for the traveling wave solutions by restricting to the case $C_2=C_3=0$ and $c>0$. This way we obtain the following family of solutions to \eqref{c6EqPDEKdV}
$$
u(x,t)=
\frac{144c e^{\sqrt{c}(C1 - x+ct)}}{e^{2\sqrt{c}(C1 - x+ct)} + 864c}.
$$

% \begin{remark}
% 	Similarly to the previous case (see Remark \ref{lodelaZK}), we may investigate the existence of explicit expressions for \eqref{primitivemodifKdV} by analyzing the roots of the polynomial $P(y)=-6y^{4}+36cy^{2}+36C_{2}y-72C_{3}$, which may lead to a complete classification of the traveling waveto the modified KdV equation. This investigation, however, lies beyond the scope of the present work.
% \end{remark}

\subsection{Case \texorpdfstring{$a(u)=\frac{u^n}{n}$}{a(u)=un/n}, with \texorpdfstring{$n\in \mathbb N$}.}
%%%%%%%%%%%%%%%%%%%%%%%%%%%%%%

For this family of KdV equations, also known as generalized KdV equation, the last Pfaffian equation \eqref{ultima1forma} to be solved takes the form:
$$
-\sqrt{-\frac{2y^{n+2}}{n(n+1)(n+2)} + cy^2 + C_2y - 2C_3}\, dz + dy
$$

Assuming $C_2=C_3=0$, we can solve this Pfaffian equation, and follow along with our procedure to obtain the following solutions to \eqref{c6EqPDEKdV}:
\begin{itemize}
	\item For $c>0$,
	$$
u(x,t)=\left(\dfrac{c n( n+1) ( n+2) }{2 \cosh^{2}\left(\frac{n}{2} \sqrt{c}  (C1 - x+ct)\right)} \right)^{\frac{1}{n}}.
$$

\item For $c<0$,
$$
u(x,t)=\left(\dfrac{c n( n+1) ( n+2) }{2 \cos^{2}\left(\frac{n}{2} \sqrt{c}  (C1 - x+ct)\right)} \right)^{\frac{1}{n}}.
$$
\end{itemize}

\subsection{Case \texorpdfstring{$a(u)=\alpha \sqrt{u}+\beta u$}{a(u)}, with \texorpdfstring{$\alpha, \beta \in \mathbb R$}.}
%%%%%%%%%%%%%%%%%%%%%%%%%%%%%%

For this choice of $a(u)$, equation \eqref{c6EqPDEKdV} becomes the Schamel--Korteweg--de Vries equation. We obtain the following family of traveling wave solutions under the assumption $C_2=C_3=0$ and $c>0$:
	$$
	u(x,t)=\left(\frac{900c e^{\frac{1}{2}\sqrt{c}(C1-x+ct)}}{240\alpha e^{\frac{1}{2}\sqrt{c}(C1-x+ct)}+e^{\sqrt{c}(C1-x+ct)}+67500\beta c+14400\alpha^{2}}\right)^2.
	$$

\subsection{Case \texorpdfstring{$a(u)=2\alpha u-\beta u^2$}{a(u)}, with \texorpdfstring{$\alpha, \beta \in \mathbb R$}.}
%%%%%%%%%%%%%%%%%%%%%%%%%%%%%%

In this case, equation \eqref{c6EqPDEKdV} becomes the Gardner equation \cite{wazwaz2007new,betchewe2013new}. The final Pfaffian equation \eqref{ultima1forma}, $\omega_1|_{\Sigma_{(C_2,C_3)}}\equiv 0$, is, in this case:
	$$
	-\frac{1}{6}\sqrt{6\beta y^{4}-24 \alpha y^3+36cy^{2}+36C_{2}y-72C_{3}}dz+dy\equiv 0.
	$$
	To effectively integrate 
	\begin{equation}\label{primitiveGardner}
		\int \frac{-6}{\sqrt{6\beta y^{4}-24 \alpha y^3+36cy^{2}+36C_{2}y-72C_{3}}} dy
	\end{equation}
	we assume, as before, $C_2=C_3=0$, and we can continue the procedure to obtain the family of solutions
	$$
u(x,t)= \frac{144c e^{\sqrt{c}(C_1-x+ct)}}{576\alpha^2 + 48e^{\sqrt{c}(C_1-x+ct)}\alpha - 864\beta c + e^{2\sqrt{c}(C_1-x+ct)}},
	$$
	for $c>0$.

\section{Conclusions}
%%%%%%%%%%%%%%%%%%%%%%%%

In this paper, we have presented a systematic, geometric approach for obtaining traveling wave solutions to the gKdV equation \eqref{c6EqPDEKdV}. The core of our methodology rests on the theory of $\mathcal{C}^{\infty}$-structures. We successfully constructed a suitable $\mathcal{C}^{\infty}$-structure for the ODE resulting from applying the traveling wave ansatz to the gKdV equation, which enabled the application of a stepwise integration algorithm. This process, which involves solving a sequence of integrable Pfaffian equations, allowed us to derive a general implicit solution for the traveling waves. 

To showcase the practical utility of our general solution, we have examined several physically significant examples of the gKdV equation. By particularizing the function $a(u)$, we were able to derive explicit solutions, including those for the modified KdV and Schamel--KdV equations, as well as for cases with general power-law nonlinearities.

It is worth noting that the search for the primitives in equation \eqref{classicKdVprimitive}, equation \eqref{primitivemodifKdV} and equation \eqref{primitiveGardner} may be further analyzed through the roots of the polynomials $P(y) = -2y^3 + cy^2 + C_2 y - 2C_3$, $P(y)=-6y^{4}+36cy^{2}+36C_{2}y-72C_{3}$ and $P(y)=6\beta y^{4}-24 \alpha y^3+36cy^{2}+36C_{2}y-72C_{3}$, respectively. This analysis would lead to the complete classification of the traveling wave solutions to those equations. In fact, this idea was successfully employed in \cite{PANCOLLANTES2025116091} to completely classify the traveling wave solutions to the modified Zakharov–Kuznetsov equation. Of course, the development of the referred analysis is beyond the scope of the present paper, and it deserves further research.

In conclusion, this work establishes the \cinf-structure formalism as an effective tool for the analysis and integration of nonlinear differential equations. The method provides a clear and systematic alternative to ad-hoc techniques and contributes to the broader understanding of integrability in nonlinear systems.
The contribution of this work, therefore, lies not in the novelty of the particular solutions, but in the successful demonstration of a general and systematic integration framework.

\section*{Declarations}

\textbf{Funding.} This research was partially supported by the grant “Operator Theory: an Interdisciplinary Approach” (ProyExcel\_00780, Junta de Andalucía), the Research Project PR2023-024 (University of Cádiz), and the research group FQM-377 (Junta de Andalucía, Spain).

\textbf{Use of generative-AI tools.} During the preparation of this work, the author used ChatGPT and Gemini exclusively for grammar and language refinement. The author reviewed and edited all content and takes full responsibility for the final version of the manuscript.

\textbf{Competing interests.} The author declares no financial or non-financial competing interests relevant to the content of this article.

\textbf{Data availability.} Not applicable.  % or specify where the data/code are stored if there is data/code

\textbf{Ethics approval / Consent to participate / Consent for publication.} Not applicable.  % only needed if your work involved human / animal subjects or identifiable data

%\section{Bibliography}
%%%%%%%%%%%%%%%%%%%%%%%%%%%%%%

% include bibliography:
\bibliographystyle{unsrturl}  
\bibliography{references.bib}

\begin{thebibliography}{10}

\bibitem{russell1844report}
J.~S. Russell.
\newblock Report on waves.
\newblock In {\em Report of the 14th meeting of the British Association for the Advancement of Science}, pages 311--390, 1844.
\newblock Held at York in September 1844.
\newblock URL: \url{http://www.ma.hw.ac.uk/~chris/Scott-Russell/SR44.pdf}.

\bibitem{korteweg1895xli}
D.~J. Korteweg and G.~de~Vries.
\newblock {XLI}. {O}n the change of form of long waves advancing in a rectangular canal, and on a new type of long stationary waves.
\newblock {\em The London, Edinburgh, and Dublin Philosophical Magazine and Journal of Science}, 39(240):422--443, 1895.
\newblock \href {https://doi.org/10.1080/14786449508620739} {\path{doi:10.1080/14786449508620739}}.

\bibitem{kdvTsutsumi}
M.~Tsutsumi, T.~Mukasa, and R.~Iino.
\newblock On the generalized {K}orteweg-de{V}ries equation.
\newblock {\em Proceedings of the Japan Academy}, 46(9):921--925, 1970.
\newblock URL: \url{https://www.jstage.jst.go.jp/article/pjab1945/46/9/46_9_921/_pdf}.

\bibitem{kdvKenig}
C.~E. Kenig, G.~Ponce, and L.~Vega.
\newblock On the (generalized) {K}orteweg-de {V}ries equation.
\newblock {\em Duke Mathematical Journal}, 59(3):585--610, 1989.
\newblock \href {https://doi.org/10.1215/S0012-7094-89-05927-9} {\path{doi:10.1215/S0012-7094-89-05927-9}}.

\bibitem{kdvChrist}
F.~M. Christ and M.~I. Weinstein.
\newblock Dispersion of small amplitude solutions of the generalized {K}orteweg-de {V}ries equation.
\newblock {\em Journal of Functional Analysis}, 100(1):87--109, 1991.
\newblock \href {https://doi.org/10.1016/0022-1236(91)90103-C} {\path{doi:10.1016/0022-1236(91)90103-C}}.

\bibitem{kalita1990modified}
B.~C. Kalita and M.~K. Kalita.
\newblock Modified {K}orteweg--de{V}ries solitons in a warm plasma with negative ions.
\newblock {\em Physics of Fluids B: Plasma Physics}, 2(3):674--676, 1990.
\newblock \href {https://doi.org/10.1063/1.859302} {\path{doi:10.1063/1.859302}}.

\bibitem{emami2011solitons}
Z.~Emami and H.~R. Pakzad.
\newblock Solitons of {KdV} and modified {KdV} in dusty plasmas with superthermal ions.
\newblock {\em Indian Journal of Physics}, 85:1643--1652, 2011.
\newblock \href {https://doi.org/10.1007/s12648-011-0178-4} {\path{doi:10.1007/s12648-011-0178-4}}.

\bibitem{tariq2022new}
K.~U. Tariq, H.~Rezazadeh, M.~Zubair, M.~S. Osman, and L.~Akinyemi.
\newblock New exact and solitary wave solutions of nonlinear {S}chamel--{K}d{V} equation.
\newblock {\em International Journal of Computer Mathematics}, 99(8):1616--1628, 2022.
\newblock \href {https://doi.org/10.1007/s40819-022-01315-3} {\path{doi:10.1007/s40819-022-01315-3}}.

\bibitem{giresunlu2017exact}
I.~B. Giresunlu, Y.~S. Ozkan, and E.~Ya{\c{s}}ar.
\newblock On the exact solutions, {L}ie symmetry analysis, and conservation laws of {S}chamel--{K}orteweg--de {V}ries equation.
\newblock {\em Mathematical Methods in the Applied Sciences}, 40(11):3927--3936, 2017.
\newblock \href {https://doi.org/10.1002/mma.4274} {\path{doi:10.1002/mma.4274}}.

\bibitem{wazwaz2007new}
A.-M. Wazwaz.
\newblock New solitons and kink solutions for the {G}ardner equation.
\newblock {\em Communications in Nonlinear Science and Numerical Simulation}, 12(8):1395--1404, 2007.
\newblock \href {https://doi.org/10.1016/j.cnsns.2005.11.007} {\path{doi:10.1016/j.cnsns.2005.11.007}}.

\bibitem{betchewe2013new}
G.~Betchewe, K.~K. Victor, B.~B. Thomas, and K.~T. Crepin.
\newblock New solutions of the {G}ardner equation: Analytical and numerical analysis of its dynamical understanding.
\newblock {\em Applied Mathematics and Computation}, 223:377--388, 2013.
\newblock \href {https://doi.org/10.1016/j.amc.2013.08.028} {\path{doi:10.1016/j.amc.2013.08.028}}.

\bibitem{inc2018investigation}
M.~Inc, A.~Yusuf, A.~I. Aliyu, and D.~Baleanu.
\newblock Investigation of the logarithmic-{KdV} equation involving {M}ittag--{L}effler type kernel with {A}tangana--{B}aleanu derivative.
\newblock {\em Physica A: Statistical Mechanics and its Applications}, 506:520--531, 2018.
\newblock \href {https://doi.org/10.1016/j.physa.2018.04.092} {\path{doi:10.1016/j.physa.2018.04.092}}.

\bibitem{pancinf-sym}
A.~J. Pan-Collantes, A.~Ruiz, C.~Muriel, and J.~L. Romero.
\newblock {$\mathcal{C}^{\infty}$}-symmetries of distributions and integrability.
\newblock {\em Journal of Differential Equations}, 348:126--153, 2023.
\newblock \href {https://doi.org/10.1016/j.jde.2022.11.051} {\path{doi:10.1016/j.jde.2022.11.051}}.

\bibitem{pancinf-struct}
A.~J. Pan-Collantes, A.~Ruiz, C.~Muriel, and J.~L. Romero.
\newblock {$\mathcal{C}^{\infty}$}-structures in the integration of involutive distributions.
\newblock {\em Physica Scripta}, 98(8):085222, 2023.
\newblock \href {https://doi.org/10.1088/1402-4896/ace403} {\path{doi:10.1088/1402-4896/ace403}}.

\bibitem{pan23integration}
A.~J. Pan-Collantes, C.~Muriel, and A.~Ruiz.
\newblock Integration of differential equations by {$\mathcal{C}^{\infty}$}-structures.
\newblock {\em Mathematics}, 11(18):3897, 2023.
\newblock \href {https://doi.org/10.3390/math11183897} {\path{doi:10.3390/math11183897}}.

\bibitem{PANCOLLANTES2025116091}
A.~J. Pan-Collantes, C.~Muriel, and A.~Ruiz.
\newblock Classification of traveling wave solutions of the modified {Z}akharov--{K}uznetsov equation.
\newblock {\em Chaos, Solitons \& Fractals}, 193:116091, 2025.
\newblock \href {https://doi.org/10.1016/j.chaos.2025.116091} {\path{doi:10.1016/j.chaos.2025.116091}}.

\bibitem{warner}
F.~W. Warner.
\newblock {\em Foundations of Differentiable Manifolds and {L}ie Groups}, volume~94 of {\em Graduate Texts in Mathematics}.
\newblock Springer-Verlag, New York, 1983.
\newblock \href {https://doi.org/10.1007/978-1-4757-1799-0} {\path{doi:10.1007/978-1-4757-1799-0}}.

\bibitem{lee2013smooth}
J.~M. Lee.
\newblock {\em Introduction to Smooth Manifolds}, volume 218 of {\em Graduate Texts in Mathematics}.
\newblock Springer, second edition, 2013.
\newblock \href {https://doi.org/10.1007/978-1-4419-9982-5} {\path{doi:10.1007/978-1-4419-9982-5}}.

\bibitem{saunders1989geometry}
D.~J. Saunders.
\newblock {\em The Geometry of Jet Bundles}, volume 142 of {\em London Mathematical Society Lecture Note Series}.
\newblock Cambridge University Press, 1989.
\newblock \href {https://doi.org/10.1017/cbo9780511526411} {\path{doi:10.1017/cbo9780511526411}}.

\bibitem{basarab}
P.~Basarab-Horwath.
\newblock Integrability by quadratures for systems of involutive vector fields.
\newblock {\em Ukrainian Mathematical Journal}, 43(10):1236--1242, 1991.
\newblock \href {https://doi.org/10.1007/BF01061807} {\path{doi:10.1007/BF01061807}}.

\bibitem{hartl1994solvable}
T.~Hartl and C.~Athorne.
\newblock Solvable structures and hidden symmetries.
\newblock {\em Journal of Physics A: Mathematical and General}, 27(10):3463--3471, 1994.
\newblock \href {https://doi.org/10.1088/0305-4470/27/10/022} {\path{doi:10.1088/0305-4470/27/10/022}}.

\bibitem{handbookfunctions}
F.~W.~J. Olver, D.~W. Lozier, R.~F. Boisvert, and C.~W. Clark, editors.
\newblock {\em {NIST} Handbook of Mathematical Functions}.
\newblock Cambridge University Press, 2010.
\newblock \href {https://doi.org/10.1080/00107514.2011.582161} {\path{doi:10.1080/00107514.2011.582161}}.

\end{thebibliography}

\end{document}